\documentclass[12pt]{amsart}

\usepackage[all]{xy}
\def\url#1{\expandafter\string\csname #1\endcsname}

\usepackage{url}

\usepackage{DLdef1}



\newcommand{\del}{\partial}
\let\ssec\subsection

\renewcommand {\ssbegin}[1]
 {\refstepcounter{subsection}
 \def \secno {\gdef \secno {}{\ssecfont \thesubsection.\hskip 2ex}%
 }%
  \begin{#1}}

\newcounter{notes}

\newcommand{\?}{\nobreak\hskip.145em\nobreak\hskip\z@skip}
\renewcommand\appendix{\par
 \setcounter{section}{0}%
 \setcounter{subsection}{0}%
 \gdef\thesection{\appendixname\kern1ex\@Alph\c@section}}

\makeatother

\begin{document}
%\graphicspath{{./figs/}}

\title[The Dzhumadildaev  brackets]{The Dzhumadildaev brackets:\\
a hidden supersymmetry of commutators\\ and the
Amitsur-Levitzki-type identities}

\author{Alexei
Lebedev${}^a$, Dimitry Leites${}^b$}

\address{${}^a$Equa
Simulation AB, 
R{\aa}sundav{\"a}gen 100, SE-169 57 Solna, Sweden; alexeylalexeyl@mail.ru\\
${}^b$Department of Mathematics, University of Stockholm, SE-106 91 Stockholm, Sweden;
mleites@math.su.se}



\begin{abstract} The Amitsur-Levitzki identity for matrices was
generalized in several directions: by Kostant for simple
finite-dimensional Lie algebras, by Kirillov (later joined by Kontsevich,
Molev, Ovsienko, and Udalova) for simple vectorial Lie algebras with
polynomial coefficients, and by Gie, Pinczon, and Ushirobira for the
orthosymplectic Lie superalgebra $\mathfrak{osp}(1|n)$.

Dzhumadildaev switched the focus of attention in these results by
considering the algebra formed by antisymmetrizors and discovered a
hidden supersymmetry of commutators.

We overview these results and their possible generalizations (open
problems).
\end{abstract}

\thanks{We are
thankful to A.~Dzhumadildaev for help and inspiring comments.}

%\date{Received May 1, 2006}

\maketitle

\markboth{Alexei Lebedev\textup{,} Dimitry Leites}{{\itshape Hidden
supersymmetry of commutators}}

\thispagestyle{empty}

\section{Introduction}\label{Sintr}
Hereafter, the ground field is $\Cee$ although several statements
are true over fields $\Kee$ of characteristic $p>2$. 

\ssec{On an experience of superizing}\label{ssExpSu} Consciously
superizing various notions and statements since 1971, people
observed that there are, usually, several ways and results of
superizations: a straightforward one (usually, not a breath-catching
one) and one or several other, often quite amazing, ways bringing up
totally new notions (examples: the Poisson and anti- brackets, the
supertrace and the queer trace on supermatrices, and the
``quasi-classical limit" of these traces, and the corresponding
superdeterminants, see \cite{LSoS} and \cite{DBS}, p. 476).

A difficulty to be able to superize something by at least one method
(to say nothing of several) usually indicates that we do not
understand, actually, even the allegedly well-understood ``nonsuper"
situation. A prime example is the integration theory on
supermanifolds which is still far from being completely constructed,
see \cite{LSoS} and \cite{Lint}. Other examples are two somewhat
related topics personified by the following two theorems:

\ssbegin{Theorem}[Cayley-Hamilton]\label{CayHa}
\begin{equation}\label{CHT}
\text{\begin{minipage}[l]{14cm} Every $n\times n$-matrix $X$
satisfies its characteristic polynomial $\det(X-\lambda 1_n)=0$.
\end{minipage}
}
\end{equation}
\end{Theorem}
Its first superization is due to Yastrebov \cite{Ya}. For various
(seemingly completely unrelated) super versions of the
Cayley-Hamilton Theorem, see \cite{KT, Del,KV,OP,GPS}.

\ssbegin{Theorem}[Amitsur--Levitzki]\label{AmiLev} Let $C$ be a
commutative and associative algebra. For any $X_{1}, \ldots ,
X_{r}\in \Mat (n;\, C)$, define antisymmetrizors $a_r$ by setting
\begin{equation}\label{a_N}
a_r(X_{1}, \ldots , X_{r}):=\sum^{}_{\sigma\in \fS_{n}} (-1)^{\sign
\sigma}X_{\sigma (1)}\ldots X_{\sigma (r)}.
\end{equation}
Then the Amitsur-Levitzki Identity (ALI) takes place:
\begin{equation} \label{6eq91}
a_r (X_{1}, \ldots ,
X_{r})=0\;\text{ for any $r\ge 2n$.}
\end{equation}
\end{Theorem}

An interesting paper \cite{GPU} was allegedly the final word
concerning superization of  ALI, but later a no less interesting
paper \cite{Sa} appeared. In this note, we also discuss
superizations of ALI; for the proof of the classical ALI with the
help of a Grassmann superalgebra, see \S\ref{ss6.6}.

\sssec{Amitsur-Levitzki type theorem for vectorial Lie algebras}
A.~A.~Kirillov formulated the following analog of the
Amitsur-Levitzki theorem, for its proof, see Preprints of Keldysh
Inst. of Applied Math. in 1980s; for a translation of one such
preprint, see \cite{KOU}; the other preprints with related results
by Kirillov, Kontsevich and Molev were never translated; Molev
reviewed them in \cite{Mo}.

\begin{Theorem}[\cite{Ki}]\label{Ki} Let $\fg$ be a simple Lie algebra of vector fields
over a field of characteristic $0$. Let
\begin{equation}\label{T_k}
A_k(x_1, \dots, x_k):=\mathop{\sum}\limits_{\sigma\in
S_k}(-1)^{\sign\sigma}\ad_{x_{\sigma(1)}} \dots \ad_{x_{\sigma(k)}}.
\end{equation}

For any $x_1, \dots, x_k\in\fg$, the identity $A_k(x_1, \dots,
x_k)\equiv 0$  holds

a) for $k\geq (n+1)^2$ if $\fg=\fvect(n)$,

b) for $k\geq n(2n+5)$ if $\fg=\fh(2n)$,

c) for $k\geq 2n^2+5n+5$ if $\fg=\fk(2n+1)$.
\end{Theorem}

\ssec{Facts that inspired us} Let $\fvect(n)$ be the Lie algebra of
vector fields (for simplicity, with polynomial coefficients).
\begin{equation}\label{F}
\begin{minipage}[c]{15cm}
%\begin{center}
\textbf{Fact}. {\sl The product of two vector fields is not a vector
field (unless is equal to $0$), but their commutator always is.}
%\end{center}
\end{minipage}
\end{equation} In \cite{D1},
Dzhumadildaev revealed a hidden supersymmetry of this well-known
Fact\eqref{F} and posed a problem natural from this super point of
view: quest for ``higher" supersymmetries on the good old Lie
algebras. Let us recall the less popular definitions and
Dzhumadildaev's construction.

Dzhumadildaev called the antisymmetrizor \eqref{a_N} of vector
fields $X_1,\ldots,X_N\in \fvect(n)$ an \emph{$N$-commutator} if
$a_N(X_1,\ldots,X_N)\in \fvect(n)$ for any $X_1,\ldots,X_N\in
\fvect(n)$ and $a_N$ does not vanish identically.  If
$a_N(X_1,\ldots,X_N)$ is an \textit{$N$-commutator}, the number
$N=N(n)$ is said to be \emph{critical}.

The $N$-commutator is \emph{subcritical} if
$A_N(X_1,\ldots,X_N):=a_N(\ad_{X_1},\ldots,\ad_{X_N})$ is
multiplication by a function for any $X_1,\ldots,X_N\in \fvect(n)$.
For example, in \cite{KOU}, it is shown that for $\fvect(1)$, the
antisymmetrizor $a_3$ acts as an operator of multiplication by a
function:
\begin{equation}\label{subcrit}
a_3(\ad_{X_1},\ldots,\ad_{X_3})(Y)=
-2\det\begin{pmatrix}x_1&x_2&x_3\\x_1'&x_2'&x_3'\\x_1''&x_2''&x_3''
\end{pmatrix}\cdot Y,
\end{equation}
where $X_i=x_i(t)\frac{d}{dt}$ for $i=1,2,3$, $f':=\frac{df}{dt}$,
and $Y=y(t)\frac{d}{dt}$.

\sssec{Problems} 1) Is the following analog to the case for $n=1$
true?
\[
\text{If $A_N(X_1,\ldots,X_N)= 0$, then $\deg
A_{N-1}(X_1,\ldots,X_{N-1})=0$.}
\]

2) The number $N(n)=2$ is always critical for any $n$; we will call
it the \emph{standard} critical number. In \cite{D1}, Dzhumadildaev
conjectured that the numbers $N(n)=(n+1)^2-3$ are also critical for
$n>1$, proved the conjecture for $n=2$ and 3, and raised a natural
problem: \textbf{List all critical numbers}. The problem is open,
except for $n=3$, where Dzhumadildaev established that $N=10$ is also
critical, and there are no more critical numbers.

Before we start considering this problem, let us discuss one more of
Dzhumadildaev's results. To present it, we need one more fact.
Although we are sure that this fact was known since at least 1960s
(for example, to I.~Kantor and/or M.~Gerstenhaber), the first
reference we know is due to Dzhumadildaev \cite{D0}:
\begin{equation}\label{F2}
\begin{minipage}[c]{14cm}
%\begin{center}
\textbf{Fact}. {\sl The antisymmetrizors form an algebra with
respect to the product defined to be
$$\renewcommand{\arraystretch}{1.4}
\begin{array}{l}(a_k\ast a_l)(X_1,\ldots,X_{k+l-1}):=\\
\mathop{\sum}\limits_{\tiny\begin{array}{l}\sigma\in S_{k+l-1}\text{
such
that}\\
\sigma (1)<\dots <\sigma (l)\;\text{ and }\\ \sigma (l+1)<\dots
<\sigma (k+l-1)\end{array}} \sign(\sigma) \, a_k(a_l(X_{\sigma
(1)}\dots X_{\sigma (l)}), X_{\sigma (l+1)}\dots X_{\sigma
(k+l-1)}).\end{array}$$}
%\end{center}
\end{minipage}
\end{equation}
More precisely, we have (\cite{D0})
\begin{equation}\label{anti}\renewcommand{\arraystretch}{1.4}
a_k\ast a_l=\begin{cases}0&\text{if $k, l$ are
even},\\
ka_{k+l-1}&\text{if $l$ is odd},\\
a_{k+l-1}&\text{if $l$ is even and $k$ is
odd.}\end{cases}\end{equation} Thus, the antisymmetrizors define a
$\Zee$-graded superring $A=\oplus A_i$, where $A_i=\Span(a_i)$, such
that $A_\ev=\mathop{\oplus}\limits_{i\equiv 1 \mod 2} A_i$ and the
product of any two odd elements of
$A_\od=\mathop{\oplus}\limits_{i\equiv 0 \mod 2} A_i$ is zero.
Clearly, $A$ can be considered as a superalgebra over any field by tensoring over $\Zee$.
\textbf{What is the meaning of the superring or superalgebra $A$?}

\ssec{Dzhumadildaev's approach to antisymmetrizors}\label{ssDzhu} In
a series of papers, Dzhumadildaev changed the emphasis of the
interpretation of the result by Amitsur and Levitzki from the search
of the identity of the least order to the description of the
superalgebra or the superring constructed from the antisymmetrizors
in the classical Lie algebras. This approach revealed a hidden
relation of the commutators with a certain universal odd
superderivation. We overview various possible generalizations of
Dzhumadildaev's result.

Let $\cF$ be an associative commutative algebra, $\cA=\End \cF$ the
associative algebra of its endomorphisms, and $\cA_L$ the Lie
algebra constructed by replacing the associative product by the
bracket. 

If\footnote{In Geometry, $\cF$ is the algebra of functions
on an $n$-dimensional manifold; it is interesting to generalize
Dzhumadildaev's approach to such cases, e.g., to the algebra of
Laurent polynomials, i.e., the algebra of functions on the torus.}
$\cF=\Kee[x_1,...,x_n]$, one can identify the elements of $\End \cF$
with differential operators. If $\End \cF$ is considered as
associative algebra, its elements satisfy no identity except
associativity. The Lie algebra $L=\Der \cF$ is a Lie subalgebra of
$(\End \cF)_L$ naturally identified with the Lie algebra $\fvect(n)$
of vector fields with polynomial coefficients.

Among numerous irreducible representations of $L$ (for their
overview, super setting including, see \cite{GLS}), there are two
``smallest" ones: in the space of functions (or, more generally,
$\lambda$-densities) and the adjoint representation.

Initially, people were interested in polynomial identities in the
adjoint representations, see Theorem~\ref{Ki}. Instead, Dzhumadildaev
considered polynomial identities in the ``smallest" representation,
which for $\fvect(n)$ is the representation in the space of
functions $\cF$. It is very interesting to generalize
Dzhumadildaev's approach on the representations in the space of
$\lambda$-densities, which is a rank 1 module over the algebra $\cF$
generated by the $\lambda$-th power of the volume element with the
following $\fvect(n)$-action (here $\lambda\in\Cee$ is fixed):
\[
X(f\vvol^\lambda)=(X(f)+f\lambda\Div(X))\vvol^\lambda \text{~~for
any $f\in \cF$ and $X\in\fvect(n)$}.
\]

It seems that this approach is more natural than the initial one for
the following reasons:

1) If one knows identities in the ``natural" representation (of the
smallest dimension or --- for infinite-dimensional algebras --- its
analog), then it is easy to construct identities in other
representations, in particular in the adjoint representation. For
example, $a_{n^2+2n-1}=0$ is identity in the space of functions
$\cF$, and since $\ad_X=r_X-l_X$, where $r_X$ and $l_X$ are right
and left actions in  $\cF$, it is easy to deduce that
$a_{n^2+2n+1}=0$ is an identity in the adjoint representation of
$\fvect(n)$.

2) If $a_N=0$ is identity, then one can ask  ``is
$a_{N-1}$ a new operation on $\fvect(n)$?"

To consider $a_{N-1}$ as a multi-operation on $\fvect(n)$  is
meaningless: $a_{N-1}$ maps $\wedge^{N-1} \fvect(n)$ to the whole
$\cA=\End \cF$, not just to $\fvect(n)$.
Dzhumadildaev suggested to consider $a_N$ on the space of
differential operators making the question ``is 
$a_{N-1}$ a new operation on $\fvect(n)$?" meaningful: in some
special cases $a_{N-1}$ maps $\wedge^{N-1} \fvect(n)$ to $\fvect(n)$
once again!

Now consider eq. \eqref{subcrit}. It means that 3-antisymmetric sum
of the adjoint derivations on $\fvect(1)$ is a multiplication
operator (not the adjoint operator). Certainly, it is an interesting
observation, but it is another topic. It has no connection with
$N$-commutators: in this setting to speak about $N$-commutator is
meaningless. Under the natural action
\[
a_3(X_1,X_2,X_3)=0 \text{~~is an identity.}
\]

Let us retell Dzhumadildaev's comments on observations due to
Kirillov, Molev, Razmuslov, Bergman, and others on identities in
$\fvect(n)$. The identities
\[
\begin{array}{ll}
a_N\equiv 0&\text{if $N\ge (n+1)^2$ for $\fvect(n)$}\\
a_N\equiv 0&\text{if $N\ge n(2n+5)$ for $\fh(2n)$}
\end{array}
\]
are not of the smallest degree. Moreover, these are ``easy"
identities. For example, for the Lie algebra $\fh(2)$ of Hamiltonian
vector fields in two indeterminates, there are two identities in
degree 7. Kirillov's identity is not minimal and it is a consequence
of these two identities. A similar situation with $\fvect(n)$.
Dzhumadildaev conjectured that the minimal identity for
representation of $\fvect(n)$ in the space of functions is of degree
$(n+1)^2-2$ whereas the degree of Kirillov's identity is
$(n+1)^2$.

\ssec{Antisymmetrizors for simple finite dimensional Lie algebras.
Exponents.}

The classical Amitsur-Levitzki theorem states that $a_{2n} = 0$ is
the minimal identity for $\fgl(n)$. For $\fo(2n+1)$ and $\fsp(2n)$,
the minimal identity is $a_{4n} = 0$; for $\fo(2n)$, the minimal
identity is $a_{4n-2} = 0$ (see \cite{AL, K1, K2}). Dzhumadildaev
formulated the following theorem (known for the serial algebras) and
gave explicit formulas for 10-antisymmetrisors in terms of the
Chevalley basis for the 7-dimensional representation of $\fg_2$.

\begin{Theorem}[\cite{D2}]\label{dzhu2} Let $A(\fg)$ be the algebra
with respect to \eqref{anti}. Then
\footnotesize{
\begin{equation}\label{A(g)}\footnotesize
\renewcommand{\arraystretch}{1.4}\arraycolsep=2pt
\begin{array}{l}
A(\fsl(n)) = \Span\{a_{2k}\mid k = 1, 2,\dots, n-1\},
\text{~~in particular, $a_{2k+1}\equiv 0$ for any $k$};\\
A(\fo(2n+1)) =\Span(\{a_{4k+1}\mid k = 1, 2,\dots, n-1\}\cup\{a_{4k+2}\mid
k = 0,1, 2,\dots,
n-1\}),\\
A(\fsp(2n)) =\Span(\{a_{4k+1}\mid k = 1, 2,\dots, n-1\}\cup\{a_{4k+2}\mid
k = 0,1, 2,\dots,
n-1\}),\\
A(\fo(2n))=\Span(\{a_{4k+1}\mid k = 1, 2,\dots, n-1\}\cup\{a_{4k+2}\mid k
= 0,1, 2,\dots,
n-2\}\cup\{a_{4n-2}\}),\\
A(\fg_2) = \Span(\{a_{2};\;\; a_{10}\}).\end{array}
\end{equation}
}\end{Theorem}

\sssbegin{Problem}\label{prob1} \emph{1)} The indices of the
antisymmetrizors are doubled exponents of the respective Lie
algebras in the cases $\fsl(n)$ and $\fg_2$, but not for $\fo$ or
$\fsp$:
\begin{equation}\footnotesize
\label{esp}\renewcommand{\arraystretch}{1.4}
\begin{tabular}{|c|l|}
\hline The Coxeter group or Lie algebra&its exponents $m_i$\cr
\hline
$A_n$ or $\fsl(n+1)$&$1, 2, 3,
\dots, n$
\cr
\begin{tabular}{l}$B_n$ or $\fo (2n+1)$ \cr
$C_n$ or $\fsp(2n)$\end{tabular} for $n\geq 2$&$1, 3, \dots,
2n-1$\cr
$D_n$ or $\fo(2n)$ &$1, 3, \dots, 2n-3; \;
n-1$\cr
$I_2^{(4)}$ or $\fg_2$ &$1, 5$\cr
$F_4$ or
$\ff_4$&$1, 5, 7, 11$\cr
$E_6$ or $\fe_6$&$1, 4, 5, 7, 8,
11$\cr
$E_7$ or $\fe_7$&$1, 5, 7, 9, 11, 13, 17$\cr
$E_8$
or $\fe_8$&$1, 7, 11, 13, 17, 19, 23, 29$\cr \hline
\end{tabular}
\end{equation}

What precisely is the relation between the indices of the
nonvanishing identically operations $a_{i}$ and the exponents?

\emph{2)} For the matrix realizations in the irreducible module
$R(\pi_1)$ of the least dimension (see the right column in table
\eqref{A}), is the following conjectural left column in table
\eqref{A} correct?
\begin{equation}\label{A}\footnotesize
\renewcommand{\arraystretch}{1.4}
\begin{tabular}{|l|l|}
\hline
$A(\ff_4) = \Span(\{a_{2},\;\; a_{10},\;\; a_{14},\;\; a_{22}\})$&$\dim R(\pi_1)=26$\\

$A(\fe_6) = \Span(\{a_{2},\;\; a_{8},\;\; a_{10},\;\; a_{14},\;\;
a_{16},\;\; a_{22}\})$&$\dim R(\pi_1)=27$\\

$A(\fe_7) = \Span(\{a_{2},\;\; a_{10},\;\; a_{14},\;\; a_{18},\;\;
a_{22},\;\; a_{26},\;\; a_{34}\})$&$\dim R(\pi_1)=56$\\

$A(\fe_8) = \Span(\{a_{2},\;\; a_{14},\;\; a_{22},\;\; a_{26},\;\;
a_{34},\;\; a_{38},\;\; a_{46},\;\; a_{58}\})$&$\dim R(\pi_1)=248$\\
\hline
\end{tabular}
\end{equation}

\emph{3)} Clearly, the algebras $A(\fg)$ may depend on the
realization of $\fg$, i.e., on the representation. And this does
happen: the algebras $A(\fsl(4))$ (corresponding to $R(\pi_1)$) and
$A(\fo(6))$
 (corresponding to $R(\pi_2)$) are different. Theorem $\ref{dzhu2}$ corresponds to
matrix realizations of the Lie algebras $\fg$ in the irreducible
module of the least (except for $\fo(6)$) dimension. \end{Problem}

\parbegin{Conjecture} For the Lie algebras with the natural matrix realization, the above
approach is reasonable. However, it seems no less reasonable to
consider Lie algebra $\fg$ embedded into their universal enveloping
algebras and look for $k$-commutators on $\fg$ inside $U(\fg)$, not
inside a particular representation. For the finite-dimensional
simple Lie algebras, only $k=2$ remains. \end{Conjecture}

The proof of Theorem $\ref{dzhu2}$ is based on the particular cases
of Lemma \ref{leosp} and \cite{K1, K2}.

\section{Superizations of Theorem $\ref{dzhu2}$}
First, let us superize the notions involved. For details of
superization, see \cite{LSoS}; we only recall here some basics. The
supermatrices are considered in the standard format. The associative
algebra $\Mat(n)$ of $n\times n$ matrices has two super analogs:
$\Mat(n|m)$ and
\begin{equation}\label{q}\footnotesize
\renewcommand{\arraystretch}{1.2}
\begin{array}{ll}
\Q(n)=&\{X\in\Mat(n|n)\mid [X, J]=0 \;\text{ for the odd invertible
operator $J$}\}=\\
&\footnotesize{\left\{\begin{pmatrix}A&B\\B&A\end{pmatrix}\mid
A,B\in \Mat(n)\right \}.}
\end{array}
\end{equation}
Accordingly, the general linear Lie algebra $\fgl(n)$ has two
superanalogs: the Lie superalgebras $\fgl(n|m)$ and $\fq(n)$
obtained from the associative superalgebras $\Mat(n|m)$ and $\Q(n)$,
respectively, by replacing the dot product by the super-bracket.

On the queer Lie superalgebra $\fq(n)$, the queer trace is defined:
\begin{equation}\label{qtr}\footnotesize
\renewcommand{\arraystretch}{1.2}
\qtr: \begin{pmatrix}A&B\\B&A\end{pmatrix}\longmapsto \tr B.
\end{equation}
The Lie superalgebra $\fsq(n)$ is the Lie subsuperalgebra of $\fq(n)$
consisting of queertraceless supermatrices.

The supermatrix $X$ is said to preserve the bilinear
form $B$ if 
\[
BX+(-1)^{p(X)p(B)}X^{st}B=0,
\] 
where the
\emph{supertransposition} $st$ describing the matrix of the dual
operator, see \cite{LSoS}, is defined as follows (in the standard
format):
\[
st: \begin{pmatrix}A&B\\C&D\end{pmatrix}\longmapsto
\begin{pmatrix}A^t&-C^t\\B^t&D^t\end{pmatrix}.
\]
Thanks to linearity, it suffices to consider only homogeneous with
respect to parity elements.

The Lie superalgebras $\fosp(n|2m)$ and $\fpe(n)$ consist of elements preserving the
nondegenerate \textbf{symmetric} bilinear form (even and odd, respectively)
the normal forms of their Gram matrices being $\diag(1_n, J_{2m})$, where
$J_{2m}=\antidiag(1_m, -1_m)$, and $J_{n|n}=\antidiag(1_n, -1_n)$,
respectively (i.e., $J_{n|n}$ coincides with $J_{2n}$ but is odd).
The same Lie superalgebras preserve \textbf{anti}symmetric nondegenerate
bilinear forms. 

For the composition $X_1 \dots X_k$ of any $k$ operators $X_1,
\dots, X_k$ (supermatrices or vector fields, or whatever) of
parities $P=(p_1,\dots,p_k)\in (\Zee/2)^k$, define its
\emph{antisymmetrizor} to be
\begin{equation}\label{sa_N}a_N(X_1,\ldots,X_N):=\mathop{\sum}\limits_{s\in S_N}
\sign(s, P) \, X_{s (1)}\dots X_{s (N)},\end{equation} where
$\sign(s, P)= \sign(s)\sign(s')$ and $s'$ is the permutation induced
by $s$ on the ordered subset of odd elements among $X_1, \dots,
X_k$. In other words, if $x_1,\dots,x_k$ are elements of a
supercommutative superalgebra whose respective parities are
$p_1+\od,\dots,p_k+\od$, then $x_{s(1)}\dots x_{s(k)}=\sign (s,P)
x_1\dots x_k$. One can express $\sign(s, P)$ in another form, more
convenient for computations. We define
\begin{equation}
\label{ssign} \sign (s,P):=\prod\limits_{1\leq i<j\leq k, \quad
s(i)>s(j)} (-1)^{p_ip_j}.
\end{equation}
Define the composition of permutations by setting
\[
s_1\circ s_2=(s_1(s_2(1)),\dots,s_1(s_2(k))).
\]
The function $\sign(s, P)$ is a 1-cocycle on $S_k$ (\cite{L}):
\[
\sign (s_1\circ s_2, P)= \sign (s_1,P) \sign (s_2,s_1(P)),~\text{
where $s_1(P)=(p_{s_1(1)},\dots,p_{s_1(k)})$.}
\]

\ssbegin{Lemma}\label{lesl} The Lie superalgebra $\fsl(m|n)$ is
closed under the $a_{2l}$ for any $m,n\geq 0$ and $l>0$. Moreover,
$a_{2l}(X_1,\dots, X_{2l})\in\fsl(m|n)$ for any $X_1,\dots,
X_{2l}\in\Mat(m|n)$. For $mn=0$, the nonvanishing identically
operations $a_{k}$ are listed in Theorem $\ref{dzhu2}$.\end{Lemma}

\begin{proof} We need to prove that $\str\,a_{2l}(X_1,\dots, X_{2l})=0$.
Let $P=(p_1,\dots,p_{2l})$ be the vector
of parities of $X_1,\dots, X_{2l}$.

For any $s\in S_{2l}$, set $s'=(s(2),\dots,s(2l),s(1))$ (i.e.,
$s'=s\circ s_0$, where $s_0=(2,\dots,2l,1)$. Then the terms in the
sum \eqref{sa_N} corresponding to $s$ and $s'$ have opposite
supertraces:
$$
\renewcommand{\arraystretch}{1.4}
\begin{array}{l}
\sign (s')\sign (s',P)\str(X_{s'(1)}\dots X_{s'(2l)})=\\
\sign (s)\sign (s_0)\sign (s,P)\sign (s_0,s(P))\str(X_{s(2)}\dots
X_{s(2l)}X_{s(1)})=\\ -(-1)^{p_1(p_2+\dots+p_{2l})}\sign (s)\sign
(s,P)\times\\ (-1)^{p_{s(1)}(p_{s(2)}+
\dots+p_{s(2l)})}\str((-1)^{p_{s(1)}(p_{s(2)}+
\dots+p_{s(2l)})}X_{s(1)}\dots X_{s(2l)})=\\ -\sign (s)\sign
(s,P)\str(X_{s(1)}\dots X_{s(2l)}). \end{array}$$ Since $2l$ --- the
order of $s_0$ --- is even, $S_{2l}$ can be represented as the
disjoint union of two sets of equal cardinalities; and the set of
elements of $S_{2l}$ can be divided in pairs of the form $(s, s\circ
s_0)$. Thus, the total supertrace of the sum \eqref{sa_N} is equal
to $0$.\end{proof}

\ssbegin{Lemma}\label{leosp} The Lie superalgebras $\fosp(m|2n)$ and
$\fpe(n)$ are closed under $a_{k}$ for $k=4l+1$ and $4l+2$ for any
$m,n,l\geq 0$.

For $\fosp(m|2n)$ and $mn=0$, the nonvanishing identically
operations $a_{k}$ are listed in Theorem $\ref{dzhu2}$; for
$\fosp(1|2n)$, we have $a_{4n}=0$ (\cite{GPU}). For $\fosp(m|2n)$
and $mn\neq 0$ but not $\fosp(1|2n)$, and for $\fpe(n)$, the $a_{k}$
for $k=4l+1$ and $4l+2$ never vanish identically.
\end{Lemma}

\begin{proof} Let $B$ be the Gram matrix of the bilinear form.
Let $X_1,\dots, X_k\in\faut(B)$ be of parities $p_1,\dots,p_k$.
%Let $B$ be the form $\fosp$ preserves.
Then $BX_i+(-1)^{p_i}X_i^{st}B=0$, and we need to show that
$$Ba_k(X_1,\dots,X_k)+(-1)^{p_1+\dots+p_k}a_k(X_1,\dots,X_k)^{st}B=0.$$

Set $s^I=(k,k-1,\dots,1)$. Then we can rewrite \eqref{sa_N} as
\[
\renewcommand{\arraystretch}{1.4}
\begin{array}{l}
a_k(X_1,\dots,X_k)=\sum\limits_{s\in S_k} \sign (s\circ s^I)\sign
(s\circ s^I,P)X_{(s\circ s^I)(1)}\dots X_{(s\circ s^I)(k)}=\\
\sum\limits_{s\in S_k} \sign (s)\sign (s^I)\sign (s,P)\sign
(s^I,s(P))X_{s(k)}\dots X_{s(1)}. \end{array}
\]
Since
\[
BX_{s(k)}\dots
X_{s(1)}=-(-1)^{p_{s(k)}}X_{s(k)}^{st}BX_{s(k-1)}\dots
X_{s(1)}=\dots=(-1)^{k+p_1+\dots+p_k}X_{s(k)}^{st}\dots
X_{s(1)}^{st}B,
\]
we have
\[
\begin{array}{l}Ba_k(X_1,\dots,X_k)=\\
\sum\limits_{s\in S_k} \sign (s)\sign (s^I)\sign (s,P)\sign
(s^I,s(P))(-1)^{k+p_1+\dots+p_k}X_{s(k)}^{st}\dots X_{s(1)}^{st}B.
\end{array}
\]
On the other hand,
\[
(X_{s(1)}\dots X_{s(k)})^{st}=\sign (s^I,s(P))X_{s(k)}^{st}\dots
X_{s(1)}^{st},
\]
so
\[
a_k(X_1,\dots,X_k)^{st}B=\sum\limits_{s\in S_k} \sign (s)\sign
(s,P)\sign (s^I,s(P))(-1)^{p_1+\dots+p_k}X_{s(k)}^{st}\dots
X_{s(1)}^{st}B.
\]
The two sums are opposite if $\sign (s^I)(-1)^k=-1$, and then
\[
Ba_k(X_1,\dots,X_k)+(-1)^{p_1+\dots+p_k}a_k(X_1,\dots,X_k)^{st}B=0.
\]
Since $\sign (s^I)=(-1)^{[k/2]}$, this is true for $k=4l+1,
4l+2$.\end{proof}

\sssbegin{Problem} What is the analog of Lemma $\ref{leosp}$ for
$\fspe(n)$?
\end{Problem}

\ssbegin{Lemma} The Lie superalgebra $\fq(n)$ is closed under
$a_{k}$ and $\fsq(n)$ is closed under $a_{2k}$ for any $n$ and $k$.
\end{Lemma}

\begin{proof} The associative algebra $\Q(n)$ is closed with respect to the
dot product; hence the result about $\fq$.

Since $\qtr(XY)=(-1)^{p(X)p(Y)}\qtr(YX)$ ($=\qtr(YX)$, since $X$ and
$Y$ should be of different parities in order to have $\qtr(XY)\neq
0$) and so the same arguments as for $\fsl$ are applicable.
\end{proof}

\ssec{Questions} What is the super analog of eq. \eqref{anti} for
the super-antisymmetrizor \eqref{sa_N}?

\section{Vectorial Lie algebras}

\ssec{$\fvect(n)$} In \cite{FF}, Feigin and Fuchs proved, among
other things, that for $n=1$, the only critical pair is the standard
one: $(1,2)$.

In \cite{D1}, Dzhumadildaev showed that for $n=2$, the complete list
of critical pairs consists only of the standard pairs $(2,2)$ and $(2,6)$. Dzhumadildaev 
gave the following explicit
expression of the 6-commutator: the 6-tuple
$(X_1,X_2,X_3,X_4,X_5,X_6)$, where
$X_i=u_{i,1}\partial_1+u_{i,2}\partial_2$ for $i=1,\ldots,6$, goes
to
\begin{equation}\label{vect6}\tiny
\begin{array}{l}
\left(\left|\begin{array}{cccccc}
u_{1,1}&u_{2,1}&u_{3,1}&u_{4,1}&u_{5,1}&u_{6,1}\\
u_{1,2}&u_{2,2}&u_{3,2}&u_{4,2}&u_{5,2}&u_{6,2}\\
\partial_2u_{1,1}&\partial_2u_{2,1}&\partial_2u_{3,1}&
\partial_2u_{4,1}&\partial_2u_{5,1}&\partial_2u_{6,1}\\
\partial_1u_{1,2}&\partial_1u_{2,2}&\partial_1u_{3,2}&
\partial_1u_{4,2}&\partial_1u_{5,2}&\partial_1u_{6,2}\\
\partial_2u_{1,2}&\partial_2u_{2,2}&\partial_2u_{3,2}&
\partial_2u_{4,2}&\partial_2u_{5,2}&\partial_2u_{6,2}\\
\partial_2^2u_{1,2}&\partial_2^2u_{2,2}&\partial_2^2u_{3,2}&
\partial_2^2u_{4,2}&\partial_2^2u_{5,2}&\partial_2^2u_{6,2}\\
\end{array}\right|+
\left|\begin{array}{cccccc}
u_{1,1}&u_{2,1}&u_{3,1}&u_{4,1}&u_{5,1}&u_{6,1}\\
u_{1,2}&u_{2,2}&u_{3,2}&u_{4,2}&u_{5,2}&u_{6,2}\\
\partial_1u_{1,1}&\partial_1u_{2,1}&\partial_1u_{3,1}&
\partial_1u_{4,1}&\partial_1u_{5,1}&\partial_1u_{6,1}\\
\partial_2u_{1,1}&\partial_2u_{2,1}&\partial_2u_{3,1}&
\partial_2u_{4,1}&\partial_2u_{5,1}&\partial_2u_{6,1}\\
\partial_2u_{1,2}&\partial_2u_{2,2}&\partial_2u_{3,2}&
\partial_2u_{4,2}&\partial_2u_{5,2}&\partial_2u_{6,2}\\
\partial_1^2u_{1,1}&\partial_1^2u_{2,1}&\partial_1^2u_{3,1}&
\partial_1^2u_{4,1}&\partial_1^2u_{5,1}&\partial_1^2u_{6,1}\\
\end{array}\right|+\right.\\
\left|\begin{array}{cccccc}
u_{1,1}&u_{2,1}&u_{3,1}&u_{4,1}&u_{5,1}&u_{6,1}\\
u_{1,2}&u_{2,2}&u_{3,2}&u_{4,2}&u_{5,2}&u_{6,2}\\
\partial_1u_{1,1}&\partial_1u_{2,1}&\partial_1u_{3,1}&
\partial_1u_{4,1}&\partial_1u_{5,1}&\partial_1u_{6,1}\\
\partial_2u_{1,1}&\partial_2u_{2,1}&\partial_2u_{3,1}&
\partial_2u_{4,1}&\partial_2u_{5,1}&\partial_2u_{6,1}\\
\partial_1u_{1,2}&\partial_1u_{2,2}&\partial_1u_{3,2}&
\partial_1u_{4,2}&\partial_1u_{5,2}&\partial_1u_{6,2}\\
\partial_2^2u_{1,2}&\partial_2^2u_{2,2}&\partial_2^2u_{3,2}&
\partial_2^2u_{4,2}&\partial_2^2u_{5,2}&\partial_2^2u_{6,2}\\
\end{array}\right|-
2\left|\begin{array}{cccccc}
u_{1,1}&u_{2,1}&u_{3,1}&u_{4,1}&u_{5,1}&u_{6,1}\\
u_{1,2}&u_{2,2}&u_{3,2}&u_{4,2}&u_{5,2}&u_{6,2}\\
\partial_2u_{1,1}&\partial_2u_{2,1}&\partial_2u_{3,1}&
\partial_2u_{4,1}&\partial_2u_{5,1}&\partial_2u_{6,1}\\
\partial_1u_{1,2}&\partial_1u_{2,2}&\partial_1u_{3,2}&
\partial_1u_{4,2}&\partial_1u_{5,2}&\partial_1u_{6,2}\\
\partial_2u_{1,2}&\partial_2u_{2,2}&\partial_2u_{3,2}&
\partial_2u_{4,2}&\partial_2u_{5,2}&\partial_2u_{6,2}\\
\partial_{12}u_{1,1}&\partial_{12}u_{2,1}&\partial_{12}u_{3,1}&
\partial_{12}u_{4,1}&\partial_{12}u_{5,1}&\partial_{12}u_{6,1}\\
\end{array}\right|-\\
2\left|\begin{array}{cccccc}
u_{1,1}&u_{2,1}&u_{3,1}&u_{4,1}&u_{5,1}&u_{6,1}\\
u_{1,2}&u_{2,2}&u_{3,2}&u_{4,2}&u_{5,2}&u_{6,2}\\
\partial_1u_{1,1}&\partial_1u_{2,1}&\partial_1u_{3,1}&
\partial_1u_{4,1}&\partial_1u_{5,1}&\partial_1u_{6,1}\\
\partial_2u_{1,1}&\partial_2u_{2,1}&\partial_2u_{3,1}&
\partial_2u_{4,1}&\partial_2u_{5,1}&\partial_2u_{6,1}\\
\partial_1u_{1,2}&\partial_1u_{2,2}&\partial_1u_{3,2}&
\partial_1u_{4,2}&\partial_1u_{5,2}&\partial_1u_{6,2}\\
\partial_{12}u_{1,1}&\partial_{12}u_{2,1}&\partial_{12}u_{3,1}&
\partial_{12}u_{4,1}&\partial_{12}u_{5,1}&\partial_{12}u_{6,1}\\
\end{array}\right|+3\left|\begin{array}{cccccc}
u_{1,1}&u_{2,1}&u_{3,1}&u_{4,1}&u_{5,1}&u_{6,1}\\
u_{1,2}&u_{2,2}&u_{3,2}&u_{4,2}&u_{5,2}&u_{6,2}\\
\partial_1u_{1,1}&\partial_1u_{2,1}&\partial_1u_{3,1}&
\partial_1u_{4,1}&\partial_1u_{5,1}&\partial_1u_{6,1}\\
\partial_1u_{1,2}&\partial_1u_{2,2}&\partial_1u_{3,2}&
\partial_1u_{4,2}&\partial_1u_{5,2}&\partial_1u_{6,2}\\
\partial_2u_{1,2}&\partial_2u_{2,2}&\partial_2u_{3,2}&
\partial_2u_{4,2}&\partial_2u_{5,2}&\partial_2u_{6,2}\\
\partial_2^2u_{1,1}&\partial_2^2u_{2,1}&\partial_2^2u_{3,1}&
\partial_2^2u_{4,1}&\partial_2^2u_{5,1}&\partial_2^2u_{6,1}\\
\end{array}\right|-\\
\left.
2\left|\begin{array}{cccccc}
u_{1,1}&u_{2,1}&u_{3,1}&u_{4,1}&u_{5,1}&u_{6,1}\\
u_{1,2}&u_{2,2}&u_{3,2}&u_{4,2}&u_{5,2}&u_{6,2}\\
\partial_1u_{1,1}&\partial_1u_{2,1}&\partial_1u_{3,1}&
\partial_1u_{4,1}&\partial_1u_{5,1}&\partial_1u_{6,1}\\
\partial_2u_{1,1}&\partial_2u_{2,1}&\partial_2u_{3,1}&
\partial_2u_{4,1}&\partial_2u_{5,1}&\partial_2u_{6,1}\\
\partial_2u_{1,2}&\partial_2u_{2,2}&\partial_2u_{3,2}&
\partial_2u_{4,2}&\partial_2u_{5,2}&\partial_2u_{6,2}\\
\partial_{12}u_{1,2}&\partial_{12}u_{2,2}&\partial_{12}u_{3,2}&
\partial_{12}u_{4,2}&\partial_{12}u_{5,2}&\partial_{12}u_{6,2}\\
\end{array}\right|
\right)\partial_1-(...)\partial_2,\\
\end{array}
\end{equation}
where the coefficient of $\partial_2$ is obtained from that of
$\partial_1$ by interchanging the subscripts $1$ and $2$ if there is
only one subscript (as in $\partial_{1}u_{5,2}$), only second subscripts 
$1$ and $2$ should be interchanged when dealing with $u_{ij}$.

\ssec{How to write the $k$-commutator for any $n$?} Let $X_1, ...,
X_k\in\fvect(n)$ with coefficients $u_{i,j}$ (i.e.,
$X_i=\mathop{\sum}\limits_{1\leq j\leq n}u_{i,j}\partial_j$). Let
$a=(a_1,\dots , a_k)$, where $a_i\in\{1, 2, \dots, n\}$; let
$(b_{ij})$ be a $k\times n$ matrix with elements in $\Zee_{\geq 0}$;
let $D((a_i), (b_{ij})$ be the determinant of the $k\times k$ matrix
whose $(i,j)$-th slot is occupied by
$\partial_1^{b_{i1}}...\partial_n^{b_{in}}u_{j,a_i}$. Considering
the $k$-commutator of the fields $X_1, ..., X_k$ as a differential
operator, its 1-st degree component is equal to
\begin{equation}\label{k-comm}\mathop{\sum}\limits_{a_1=1}^n\dots\mathop{\sum}\limits_{a_k=1}^n
\mathop{\sum}\limits_{s_1=2}^k \mathop{\sum}\limits_{s_2=3}^k
\dots\mathop{\sum}\limits_{s_{k-1}=k}^k
D\left((a_i),\mathop{\sum}\limits_{1\leq i\leq k-1}
E^{s_i,a_i}\right)\del_{a_k},
\end{equation} where the $E^{i,j}$ are matrix units.

Accordingly, if the $k$-commutator is a first order operator, then
\eqref{k-comm} is its expression. Unfortunately, this expression is
not user-friendly: first, it is longish ($n^k\times(k-1)!$ summands)
which even for $n=2, k=6$ is $>7000$), second, it is very redundant:
some of the summands vanish, some are equal to each other, some are
equal in absolute value but are of different signs (so there are
just 14 distinct types of summands for $n=2$, not $>7000$).

In \cite{D3}, Dzhumadildaev showed that for $n=3$, in addition to
the standard pairs $(3,2)$ and $(3,13)$, there is exactly one more
critical pair, $(3,10)$.

\ssec{The other series of simple vectorial Lie algebras with
polynomial coefficients.} It is equally natural to list all critical
pairs for the other types of simple vectorial Lie algebras. For
these Lie algebras, only the pairs $(n,2)$ will be called
\emph{standard}.

For the Lie algebras $\fsvect(n)$ of divergence-free vector fields,
Dzhumadildaev proved \cite{D1, D3} that the only nonstandard
critical pairs are $(2, 5)$ and $(3,10)$ (for $n=2$ and 3,
respectively). Since $\fsvect(2)\simeq\fh(2)$ the result for this
Lie algebra might be pertaining to the Hamiltonian series, rather
than to the divergence-free one.

For the Lie algebras $\fh(2n)$ of Hamiltonian vector fields,
Dzhumadildaev proved \cite{D1} that the only nonstandard critical
pair for $n=1$ is $(2, 5)$. In terms of generating functions in $p$
and $q$, the 5-commutator is proportional to the following beautiful
map
\begin{equation}\label{h5}
\footnotesize
\renewcommand{\arraystretch}{1.4}
 (f_1, f_2,f_3,f_4,f_5)\mapsto\det\begin{pmatrix}
 \partial_q(f_1)& \partial_q(f_2)&\partial_q(f_3)&\partial_q(f_4)&\partial_q(f_5)\\
 \partial_p(f_1)& \partial_p(f_2)&\partial_p(f_3)&\partial_p(f_4)&\partial_p(f_5)\\
 \partial_p^2(f_1)& \partial_p^2(f_2)&\partial_p^2(f_3)&\partial_p^2(f_4)&\partial_p^2(f_5)\\
 \partial_q^2(f_1)& \partial_q^2(f_2)&\partial_q^2(f_3)&\partial_q^2(f_4)&\partial_q^2(f_5)\\
 \partial_p\partial_q(f_1)& \partial_p\partial_q(f_2)&\partial_p\partial_q(f_3)&
 \partial_p\partial_q(f_4)&\partial_p\partial_q(f_5)\end{pmatrix}.
\end{equation}
%where for any $f$ and any operator $D$, we set $Df:=(Df_1,
%Df_2,Df_3,Df_4,Df_5)$.

\sssbegin{Problem} What are the $N$-commutators for the Lie algebra
of contact vector fields $\fk(2n+1)$?\end{Problem}

Dzhumadildaev's guiding idea is
very simple and brought to the title of \cite{D3}: it is a certain
odd derivation of a certain superalgebra associated with the
problem, which is in the heart of this matter, see the next Section.

\section{The universal odd derivation and $N$-commutators (after \cite{D3})}
Let $L$ be a Lie (super)algebra, $U(L)$ its enveloping algebra,
$\Pi$ the change of parity functor. Take the associative
supercommutative superalgebra $K=S^{\bcdot}(\Pi(L))$; in particular,
if $L$ is purely even, then $K$ is a Grassmann superalgebra. In $L$,
select an arbitrary basis $B$ and set
$$D=\mathop{\sum}\limits_{b\in B} \Pi(b)\otimes b \in K\otimes L\subset K\otimes U(L).$$

\ssbegin{Lemma} The $N$-commutator on $L$ yields an element of $L$,
if and only if $D^N\in K\otimes L$. The $N$-commutator does not
vanish identically if and only if $D^N\neq 0$.

In particular, for $L\subset\fvect(n)=\fder(\Kee[x])$, we clearly
have
$$D\in K\otimes L\subset \fvect(n|n)=\fder(\Kee[x, \Pi(x)]).$$ The
$N$-commutator on $L$ yields an element of $L$ if and only if
$D^N\in \fvect(n|n)=\fder(\Kee[x, \Pi(x)])$.
\end{Lemma}

\textbf{Comment}. For superspaces, the following modification of
Fact \eqref{F} takes place:
\begin{equation}\label{F3}
\begin{minipage}[c]{14cm}
%\begin{center}
\textbf{Fact}. {\sl The product of two nonproportional odd vector
fields is usually not a vector field, but the square of any odd
field is always a vector field.}
%\end{center}
\end{minipage}
\end{equation}
Fact \eqref{F} is, therefore, a corollary of Fact \eqref{F3} for
$N=2$. \textbf{This is the hidden supersymmetry of the
anticommutator mentioned in the title of this paper}.

\sssbegin{Conjecture} We only considered Lie algebras of vector
fields with polynomial coefficients. We conjecture that the answer
will be same for any type of coefficients (at least, if polynomials
are dense in the space of coefficients). \end{Conjecture}

\ssec{Discussion and setting of the problem.} As noted in Introduction, attempts to
superize a problem or a notion usually reveal two roads: a straightforward
one (not of much interest) and a totally unexpected one. The problem Dzhumadildaev 
posed (\textbf{describe all critical
pairs for (simple) Lie algebras of vector fields}) is the one which we do not know 
how to superize. In particular,
what is the answer for any of the simple Lie superalgebras of vector fields (with
polynomial coefficients to begin with)?

Recall steps of Dzhumadildaev's proof.
Let $l$ be the length function on $\Diff_n$ defined by
Dzhumadildaev, namely:
$$
l((\eta_{i_1,\;\alpha_1}\del^{\beta_1})\dots
(\eta_{i_k,\;\alpha_k}\del^{\beta_k})):=k.
$$
Let us extend $l$ to a grading (Dzhumadildaev's definition is
slightly different but equivalent). Note that the possibility of
such extension is a little less evident than in the case of
$\mathcal{L}_n$ because the elements
$\eta_{i_j,\alpha_j}\del^\beta_j$ do not supercommute.

Let $X_1,\dots, X_k$ be some abstract vector fields (considered as
variables here) of $n$ indeterminates. Define the following map $F$
from $\Diff_n^{[k]}$ to the algebra of differential operators (of
arbitrary degree) in $n$ indeterminates:
$$
F((\eta_{i_1,\;\alpha_1}\del^\beta_1)\dots
(\eta_{i_k,\;\alpha_k}\del^\beta_k))=\sum\limits_{s\in S_k}
(-1)^{\sign (s)}
((\del^{\alpha_1}X_{s(1),\;i_1})\del^{\beta_1})\dots
((\del^{\alpha_k}X_{s(k),\;i_k})\del^{\beta_k}).
$$
Here $(\del^{\alpha_j}X_{s(j),\;i_j})$ is a function,
$(\del^{\alpha_j}X_{s(j),\;i_j})\del^{\beta_j}$ is a differential
operator (possibly of zero degree, if $\beta_j=0$), and the whole
term is the composition of differential operators.

\ssbegin{Statement} The map $F$ is faithful.\end{Statement}

The idea of a proof: the map preserves commutation relations. Note
that

a) $F(D^k)$ is just the $k$-commutator of the $X_j$ (considered as a
differential operator of arbitrary degree);

b) the map $F$ preserves the degree of the differential operator
(for generic $X_i$).

So the $k$-commutator is of degree $1$ for any $X_j$ if and only if
$\deg D^k=1$.

\section{Appendix: A proof of the classical Amitsur--Levitzki identity}
\label{ss6.6}

Let $A$ be a~supercommutative superalgebra and $X\in \Mat(n|0;\,
A)_\od$. It is clear that $X^{r}=0$ for any $r>n^{2}$. It turns out
that $r$ can be  considerably diminished.

\ssbegin{Proposition}\label{prop} $X^{2n}=0$ for any $X\in \Mat(n|0;\,
A)_\od$. \end{Proposition}

First of all, let us discuss what does this identity mean from the
``ordinary", i.e., nonsuper, algebra point of view. Let $C$ be
commutative algebra and $X_{1}, \ldots, X_{r}\in \Mat (n;\, C)$. Set
$A=C[\xi _{1}, \ldots, \xi _{r}]$, where the $\xi _{i}$ are odd and
let $X:=\sum \xi_{i}X_{i}\in \Mat (n|0;\, A)_\od$. Clearly,
\begin{equation}
\label{6eq89} X^{r}= a_r (X_{1}, \ldots, X_{r})\xi_{1}\ldots
\xi_{r},\text{~~ where $a_r(X_{1}, \ldots ,
X_{r})=\sum^{}_{\sigma\in \fS_{n}} (-1)^{\sgn \sigma}X_{\sigma
(1)}\ldots X_{\sigma (r)}$}.
\end{equation}
Hence, Proposition \ref{prop} implies the Amitsur--Levitzki identity
\eqref{6eq91}.

\sssbegin{Exercise} The Amitsur--Levitzki identity implies Proposition \ref{prop}.
\end{Exercise}

\begin{proof}[Proof\nopoint] of Proposition \ref{prop}.\footnote{This proof is due to
J.~Bernstein, 1975. At about the same time V.~Drinfeld also noticed
the equivalence proved here.} Set $Y=X^{2}$. The elements of $Y$
belong to the commutative algebra $A_{\ev}$, and therefore, we may
consider the characteristic polynomial $P(\lambda)=\det (\lambda 1_n
-Y)$ with coefficients in $A_\ev$. Let us prove that
$P(\lambda)=\lambda^{n}$. Since the Cayley--Hamilton theorem implies
$P(Y)=0$, we have $Y^{n}=0$, i.e., $X^{2n}=0$. We will prove that
$P(\lambda)=\lambda^{n}$ by three different methods.

1) If $\text{char }k=0$, then the coefficients of $P(\lambda)$ can
be expressed in terms of $\tr Y^{r}$ for $r=1, 2, \ldots $.
Therefore, it suffices to verify that $\tr Y^{r}=0$. Indeed,
\begin{equation}
\label{6eq92} \tr Y^{r}=\str X^{2r}=\str X\cdot X^{2r-1}=-\str
X^{2r-1}\cdot X=-\tr Y^{r}.
\end{equation}
Hence, $\tr Y^{r}=0$ for $r=1, 2, \ldots $.

2) Let us show that $P(\lambda )^{2}=\lambda ^{2n}$. If $2$ is
invertible in $A$, we see that $P(\lambda )=\lambda ^{n}$. We have
to show that $\det ^{2}(1_n-\lambda X^{2})=1$. This follows from
a~more general statement.

\sssbegin{Lemma} Let $U\in \Mat (p\times q;\, A)$ and $V\in \Mat (q\times p;\, A)$
be matrices whose entries are odd elements of $A$. Then
\begin{equation}
\label{6eq93} \det (1_{p}-UV)=\det (1_{q}-VU)^{-1}.
\end{equation}
 \end{Lemma}

\begin{proof} Let $Z=\mat{1_p&U\cr V&1_q} \in
\GL(p|q;\, A)$ and
$\mat{A&B\\ C&D}^\Pi:=\mat{D&C\\
B&A}$. From \cite{LSoS} we know that $\Ber Z^\Pi=(\Ber Z)^{-1}$, so
$\Ber Z^\Pi=\det (1_{q}-VU)$ because  $\Ber Z=\det (1_{p}-UV)$.
\end{proof}

 3) Let $Z=\mat{1_n&\lambda X\cr \lambda X&1_n}
\subset \GQ(n;\, A[\lambda ])$. From \cite{LSoS} we know that $\Ber
Z=1$. But $\Ber Z=\det (1_n-\lambda ^{2}X^{2})$, hence, $\det
(1_n-\lambda ^{2}Y)=1$, and we have $\det (1_n-\lambda Y)=1$. Thus,
$\det (\lambda1_n -Y)=\lambda ^{n}$.
\end{proof}

\ssec{How to superize the Cayley--Hamilton theorem?}\label{ss6.6.1}
The degree of the polynomial equation a~given $n\times n$ matrix
satisfies given by the Amitsur--Levitzki identity can be diminished
even more (\textbf{Cayley--Hamilton theorem}, see \eqref{CHT}).

\sssbegin{Conjecture}[ \cite{GPS}] The analog of the Cayley--Hamilton theorem for
supermatrices was unknown, except for small values of $n$ (equal to
$2$ or $1|1$), until recently. Now we have a~conjectural formula
suggested by the study of quantum algebras and passage to the
appropriate ``super" limit.\end{Conjecture}

%\label{lastpage}


\begin{thebibliography}{2000}

\bibitem[AL]{AL}
Amitsur A., Levitzki J., Minimal identities for algebras. Proc.
Amer. Math. Soc., 1, (1950) 449--463


\bibitem[Del]{Del}
Deligne P., CatŽegories tensorielles. Moscow Math. J. 2 (2) (2002)
227-248.


\bibitem[DBS]{DBS}
Duplij S., Bagger J., Siegel W. (eds.)  \textit{Concise Encyclopedia
of~Supersymmetry and Noncommutative Structures in Mathematics and
Physics}, Kluwer, Dordreht, 2003

\bibitem[D0]{D0}
Dzhumadil'daev A.,
Integral and mod $p$-cohomologies of the Lie algebra $W_1$. Funct.
Anal. Pril. 22 (3)(1988), 68-70 (in Russian) English translation in:
Funct. Anal. Appl. 22 (1988), No.3, 226--228



\bibitem[D1]{D1}
Dzhumadil'daev A., $N$-commutators, Comment. Math. Helv. \textbf{79}
(2004), no.~3, 516--553; \texttt{arXiv:math/0203036}



\bibitem[D2]{D2}
Dzhumadil'daev A., $N$-commutators for simple Lie algebras. In: P.
Kielanowski, A. Odzijewicz, M. Schlichenmaier, Th. Voronov (Eds.)
\emph{XXVI International Workshop on Geometrical Methods in Physics:
Bialoweza, Poland 1--7 July 2007}. Amer. Inst. of Physics, (2007)
159--168

\bibitem[D3]{D3}
Dzhumadil'daev A., $10$-commutators, $13$-commutators, and odd
derivations. Journal of Nonlinear Mathematical Physics. Volume 15,
Number 1 (2008), 87-103; \texttt{arXiv:math-ph/0603054}


\bibitem[D4]{D4}
Dzhumadil'daev A., $2p$-commutator on differential operators of
order $p$. Letters in Math. Phys.  V. 104, Issue 7, (2014) 849--869;
\texttt{arXiv:1401.1730}

\bibitem[FF]{FF}
Feigin B.L., Fuchs D.B., Invariant skew-symmetric differential
operators on the line and Verma modules over the Virasoro algebra,
Funct. Anal. Appl. 16 (1982), 114--126.

\bibitem[GLS]{GLS}
Grozman P., Leites D., Shchepochkina I., Invariant differential
operators on supermanifolds and The Standard Model. In: Olshanetsky
M., Vainshtein A., (eds.) \emph{M.~Marinov memorial volume}, World
Sci., 2002, 508--555. \texttt{arXiv:math.RT/0202193}

\bibitem[GPU]{GPU}
Gie P.-A., Pinczon G., Ushirobira R., Back to the Amitsur-Levitzki
theorem: a super version for the orthosymplectic Lie superalgebra
$\fosp(1, 2n)$. Lett. Mathem. Physics, V. 66, no.~1--2, 2003,
141--155; \texttt{arXiv:math/0309418v2}

\bibitem[GPS]{GPS}
Gurevich D., Pyatov P., Saponov P., Cayley-Hamilton theorem for
quantum matrix algebras of $GL(m|n)$ type. Algebra i Analiz, 17,
no.~1 (2005) 160-181 (in Russian). \texttt{arXiv:math.QA/0412192}

%\bibitem[IMSZ]{IMSZ}
%Iuliu-Lazaroiu C., McNamee D., S\"amann Ch., Zejak A., Strong
%Homotopy Lie algebras, generalized Nahm equations and multiple
%$M2$-branes; (2009) \texttt{arXiv:hep-th/0901.3905}

\bibitem[Kag1]{Kag1}
Kagarmanov, A.A.
The standard Lie polynomial of degree 8 on the Lie algebra $W\sb 2$.
(English. Russian original) Mosc. Univ. Math. Bull. 45, no.6
(1990), 34--35

\bibitem[Kag2]{Kag2}
Kagarmanov, A. A. On the nonexistence of a polynomial reproducing
the function ring from the Lie algebra of vector fields over fields
of positive characteristic. Russ. Math. Surv. 53, no.~1, (1998) 212--213


\bibitem[KaR]{KaR}
Kagarmanov, A.A., Razmyslov, Yu.P., Existence of central polynomials
for adjoint representations of simple Lie algebras. Fundanientalnaya
i prikladnaya matematika. (Russian. English summary) Fundam. Prikl.
Mat., vol. 5 (1999), no.~4, 1015--1025 (in Russian)

\bibitem[KT]{KT}
Kantor I., Trishin I., On the Cayley-Hamilton equation in the
supercase. Comm. in Algebra, 27 (1999) 233--259

\bibitem[Ki]{Ki}
Kirillov A., On identities in the Lie algebras of vector fields.
Vestnik Mosk. universiteta. Ser. 1 Matematika i mehanika. (1989)
no.~2, 11--17 (in Russian)


\bibitem[KOU]{KOU}
Kirillov A., Ovsienko V., Udalova O., Identities in the Lie algebras
of vector fields on the real line. Selecta Mathematica Sovietica,
v.10 no. 1 (1991), 7--17


\bibitem[KV]{KV}
Khudaverdyan H., Voronov Th., Berezinians, exterior powers and
recurrent sequences, Letters in Math. Physics, 74(2), (2005), 201--228


\bibitem[K1]{K1}
Kostant B., A theorem of Frobenius, a theorem of Amitsur-Levitzki
and cohomology theory. J. Math. Mech., 7 (1958) 237--264

\bibitem[K2]{K2}
Kostant B., A Lie algebra generalization of the Amitsur-Levitzki
theorem, Adv. in Math. 40 (1981), 155--175.

%\bibitem[LM]{LM}
%Lada T., Markl M. Strongly homotopy Lie algebras. Commun. in Algebra
%23 (1995) 2147--2161; \texttt{arXiv:hep-th/9406095}

\bibitem[L]{L}
Leites D., Cohomology of Lie superalgebras, Functional Anal. Appl.,
9 (1975), 340-341

\bibitem[Lint]{Lint}
Leites D., On unconventional integration on supermanifolds and cross
ratio on classical superspaces. In: E.~Ivanov, S.~Krivonos,
J.~Lukierski, A.~Pashnev (eds.) \textit{Proceedings of the
International Workshop ``Supersymmetries and Quantum Symmetries",
September 21-25, 2001, Karpacz, Poland}. JINR, Dubna, (2002),
251--262; \texttt{arXiv:math.RT/0202194}

\bibitem[LSoS]{LSoS}
Leites D. (ed.) \emph{Seminar on supersymmetries}, (J.~Bernstein,
D.~Leites, V.~Molotkov, V.~Shander), MCCME, Moscow, (2011), 410 pp (in
Russian; an English version in preparation is available for
perusal)


\bibitem[Mo]{Mo}
Molev A., On the algebraic structure of the Lie algebra of vector
fields on the line. Mathematics of the USSR-Sbornik, (1989), 62:1,
83--94


\bibitem[OP]{OP}
Ogievetsky O., Pyatov P., Orthogonal and
Symplectic Quantum Matrix Algebras and Cayley-Hamilton Theorem for
them. \texttt{arXiv:math/0511618}



\bibitem[Sa]{Sa}
Samoilov L.M. An analog of the Amitsur-Levitzki theorem for matrix
superalgebras. Siberian Mathematical Journal, (2010), 51:3, 491--495


\bibitem[Ya]{Ya}
Yastrebov A.V. Cramer's and Cayley-Hamilton theorems for matrices
over a superalgebra. In: A.~Onishchik (ed.) \textit{Questions of
group theory and homology algebra}, Yaroslavl Univ. Press,
Yaroslavl, (1988), 130--141. (in Russian)

\end{thebibliography}
\end{document}